# Relaxation for highly discontinuous, possibly unbounded, integral functionals

T. Bertin[*]    G. Treu[†]

October 20, 2025

## Abstract

We consider the functional

$$F(u) = \int_\Omega f(\nabla u)\, dx \qquad u \in \varphi + W_0^{1,1}(\Omega)$$

where $\Omega$ is a Lipschitz bounded open set of $\mathbb{R}^N$, $f : \mathbb{R}^N \to \mathbb{R} \cup \{+\infty\}$ is a superlinear Borel function, $\varphi \in W^{1,\infty}(\Omega)$.

We prove that, if $f$ is superlinear and satisfies very weak assumptions, then the Lavrentiev phenomenon does not occur. We underline that our assumptions include the case of non continuous, non convex, and unbounded Lagrangians.

## 1 Introduction

The Lavrentiev phenomenon, first discovered by M. Lavrentiev in 1926 [36], represents a fundamental problem in the Calculus of Variations. It occurs when the infimum of a variational functional taken over a smaller, smoother class of admissible functions, for instance Lipschitz or $\mathcal{C}^1$ functions, is strictly greater than the infimum taken over a larger, weaker class, such as Sobolev functions. In other words, the expected approximation of minimizers by smooth functions fails, revealing a gap between the two infima and highlighting situations where regularity and approximation properties break down. Several examples of its occurrence have been exhibited in various cases of the Calculus of Variations also for very regular Lagrangians [1, 3, 4, 25, 31, 32, 37, 40, 49].

[*]Dipartimento di Matematica "Tullio Levi-Civita", Università di Padova, via Trieste 63 I-35121, Padova, bertin@math.unipd.it

[†]Dipartimento di Matematica "Tullio Levi-Civita", Università di Padova, via Trieste 63 I-35121, Padova, giulia.treu@unipd.it

The problem of determining sufficient conditions for the non occurrence of the Lavrentiev gap has been studied in various settings, including problems with non-standard growth conditions, constraints, or degeneracies in the integrand. We cite here some of the results in this direction, without claiming to be exhaustive [2, 5, 6, 8, 9, 10, 12, 13, 14, 20, 27, 28, 30, 31, 35, 38, 39, 40, 42, 49]

Understanding the mechanisms that give rise to or prevent the Lavrentiev phenomenon not only deepens the theoretical framework of the calculus of variations but also has implications for numerical analysis, materials science, and the theory of partial differential equations. Recently, its implications in AI have been highlighted. We underline that, while most of the results cited above deal with Lagrangians that are convex in the gradient variable, in the framework of AI non convex problems play a crucial role.

In this work, we consider the functional

$$F(u) = \int_\Omega f(\nabla u)\, dx \qquad u \in \varphi + W_0^{1,1}(\Omega) \tag{1.1}$$

where $\Omega \subset \mathbb{R}^N$ is an open bounded Lipschitz domain, $f : \mathbb{R}^N \to \mathbb{R} \cup \{+\infty\}$ is a Borel function, $\varphi$ is assumed to be Lipschitz.
In Theorem 4.3, we prove that the functional (1.1) does not exhibit the Lavrentiev phenomenon as soon as we assume $f$ is superlinear and satisfies a very weak assumption, see (HF') (ii) and (iii). We remark that $f$ is not necessarily continuous, convex, or bounded, as it has been underlined in Example 2.1.

The main difficulty in the proof of this result consists in the fact that we cannot use approximations of the function $u$ by means of convolutions. Instead, we must explicitly construct a suitable approximating sequence that has the property that each function in the sequence has a gradient taking values where $f$ and its lower semicontinuous convex envelope $f^{**}$ coincide. This is achieved in Theorem 3.1 where, refining a construction by Cellina [22] (see also [34]) later developed in [47, 45, 46], we show that, for any almost everywhere differentiable Sobolev function $u$ and any $\varepsilon > 0$, we can define a suitable function $v$, coinciding with $u$ on

$\partial\Omega$, such that
$$F(v) \leq F^{**}(u) + \varepsilon$$
where
$$F^{**}(v) = \int_\Omega f^{**}(\nabla v)\, dx.$$

The proof of the non occurrence of the Lavrentiev phenomenon for the functional $F$ then works as follows. First of all, we exploit some recent results [10, 12] that state the non occurrence of the gap for autonomous convex problems, applying them to the functional $F^{**}$. Then, denoting by $u_n$ a sequence of Lipschitz functions that is minimizing for $F^{**}$, we apply Theorem 3.1 to each element $u_n$, obtaining a new sequence of Lipschitz $v_n$ that is minimizing for the functional $F$.

## 2 Notation, preliminaries and assumptions

We denote by $|\cdot|$ the Euclidean norm in $\mathbb{R}^N$ and by $\|\cdot\|_\infty$ the norm in $L^\infty$; by $e_i$ the $i$-th vector of the canonical basis in $\mathbb{R}^N$ and by $a \cdot b$ the scalar product between two vectors $a, b \in \mathbb{R}^N$. In the whole paper $f$ will denote a Borel function $f : \mathbb{R}^N \to \mathbb{R} \cup \{+\infty\}$.

We will use the lower semicontinuous convex envelope of $f$ i.e., the function $f^{**} : \mathbb{R}^n \to \mathbb{R} \cup \{+\infty\}$. Recalling ([23],[29]) can be represented as

$$f^{**}(\xi) = \inf\left\{ \sum_{i=1}^{N+1} \alpha_i f(\xi_i) : \sum_{i=1}^{N+1} \alpha_i \xi_i = \xi,\ \alpha_i \geq 0,\ \sum_{i=1}^{N+1} \alpha_i = 1 \right\}$$

or also as
$$f^{**}(\xi) = \sup\left\{ l(\xi) : l \text{ affine},\ l \leq f \right\}.$$

We assume that $f$ satisfy the following properties.

(HF) (i) There exist $c_1 > 0$, $c_2 \in \mathbb{R}$ such that
$$f(\zeta) \geq c_1 |\zeta| + c_2 \qquad \text{for every } \zeta \in \mathbb{R}^N$$

(ii) For every $\xi \in \mathbb{R}^N$ either $f^{**}(\xi) = f(\xi)$ or there exist $k \in \{1, \ldots, N\}$, $\{\xi_i\}_{i=1\ldots,k+1} \subset \mathbb{R}^N$ and $\{\alpha_i\}_{i=1\ldots,k+1} \subset \mathbb{R}$ such that $\alpha_i$ is strictly positive for every $i = 1, \ldots, k+1$, $\sum_{i=1}^{k+1} \alpha_i = 1$, $\sum_{i=1}^{k+1} \alpha_i \xi_i = \xi$ and
$$\sum_{i=1}^{k+1} \alpha_i f(\xi_i) = f^{**}(\xi).$$

Moreover we also assume that $\dim \mathrm{Span}(\xi_1, ..., \xi_{k+1}) = k$.

(iii) For every $R > 0$ there exist $N+1$ vectors $\zeta_j \in \mathbb{R}^N$ such that $\overline{B}(0, R) \subset \overline{\mathrm{co}}(\cup_{j=1}^{N+1}\{\zeta_j\})$ and

$$f(\zeta_j) = f^{**}(\zeta_j) < +\infty \qquad \text{for every } j = 1, ..., N+1.$$

The first of the following examples will show that the assumption (HF) is very general and includes functions that assume the value $+\infty$ on very large sets, possibly a.e. on $\mathbb{R}^N$. The second one will exhibit a function that does not satisfy (HF) (ii).

**Example 2.1.** *Condition (HF) (ii) is satisfied by any function $f : \mathbb{R}^N \to \mathbb{R}$ such that for any face $\mathcal{F}$ of the epigraph of $f^{**}$ it holds*

$$\mathcal{F} = \mathrm{co}\{(\eta, f^{**}(\eta)) : f^{**}(\eta) = f(\eta) \text{ and } (\eta, f^{**}(\eta)) \in \mathcal{F}\}$$

*These class includes, for instance, functions that could be equal to $+\infty$ almost everywhere. The function defined by*

$$f(\eta) = \begin{cases} |\eta|^2, & \text{if } \eta \in \mathbb{Z}^N, \\ |\eta|^2 + g(\eta), & \text{otherwise.} \end{cases} \qquad (2.1)$$

*where $g : \mathbb{R}^n \to \mathbb{R} \cup \{+\infty\}$, $g(\eta) \geq 0$, satisfies assumption (HF) (ii).*

**Remark 2.2.** *We underline the fact that, in the case $N = 2$, the function defined in (2.1) assumes the value $+\infty$ on every line passing through the origin with an irrational slope.*

**Example 2.3.** *We consider the function defined by*

$$f(\eta) = g(|\eta|),$$

*where*

$$g(t) = \begin{cases} |t^2 - 1|, & \text{if } t \neq 1, \\ 1, & \text{if } t = 1, \end{cases}$$

*then $f$ does not satisfy condition (HF) (ii). In fact $f^{**}(\eta) = 0$ for every $\eta \in \overline{B(0, 1)}$ while $f(\eta) > 0$ for every $\eta \in \mathbb{R}^N$.*

In the proofs of our results we will use some consequences of assumption (HF).

**Lemma 2.4.** *Let $f : \mathbb{R}^n \to \mathbb{R} \cup \{+\infty\}$ satisfy assumption* (HF). *Then the following properties hold.*

(a) *The effective domain of $f^{**}$ is $\mathbb{R}^N$, i.e., $f^{**} : \mathbb{R}^N \to \mathbb{R}$;*

(b) *If $\xi$, $\xi_i$, $\alpha_i$ are as in* (HF) (ii), *then $f(\xi_i) = f^{**}(\xi)$ and there exist $a \in \mathbb{R}^n$ such that $a \in \partial f^{**}(\xi) \cap (\cap_{i=1,\ldots,k+1} \partial f^{**}(\xi_i))$.*

To conclude this section, we list some further notations that we will use in the paper. Since we will need to use Lebesgue measures in spaces of different dimensions, when dealing with the measure in $\mathbb{R}^N$, we write $\lambda$; whereas when dealing with the measure in a subspace of dimension $k$, we use $\lambda_k$. We use the standard notations for the Sobolev spaces $W^{1,1}(\Omega)$, $W^{1,\infty}(\Omega)$, $W_0^{1,1}(\Omega)$ and $W_0^{1,\infty}(\Omega)$. In particular, we denote by $\rightharpoonup$ the weak convergence in $W^{1,1}(\Omega)$, by $\stackrel{*}{\rightharpoonup}$ the weak* convergence in $W^{1,\infty}(\Omega)$, and by $\|\cdot\|_{W^{1,1}}$ and $\|\cdot\|_{W^{1,\infty}}$ respectively the norms in $W^{1,1}(\Omega)$ and in $W^{1,\infty}(\Omega)$.

## 3 An approximation result

In this section, we prove a theorem that constitutes the first step in identifying sufficient conditions to prevent the Lavrentiev phenomenon for integral functionals with non-convex and discontinuous Lagrangians. Given a function $u \in W_\varphi^{1,1}$, assumed also to be differentiable almost everywhere in $\Omega$, we will construct a function $v$ whose gradients take values where the functions $f$ and $f^{**}$ coincide. The construction is carried out so that the value of $F(v)$ is not far from $F^{**}(u)$.

**Theorem 3.1.** *Let be $\Omega$ a open bounded Lipschitz subset of $\mathbb{R}^N$ and assume that $f : \mathbb{R}^N \to \mathbb{R} \cup \{+\infty\}$ satisfies assumption* (HF). *Let $u \in \varphi + W_0^{1,1}(\Omega)$ be differentiable almost everywhere.*
*Then for every $\varepsilon > 0$ there exists $v \in \varphi + W_0^{1,1}(\Omega)$, differentiable almost everywhere, such that*
$$\int_\Omega f(\nabla v)\,dx \leq \int_\Omega f^{**}(\nabla u)\,dx + \varepsilon\,.$$

**Proof.** The construction of the function $v$ proceeds in several steps. We start by defining a piecewise affine function that will be used to modify the function $u$ in a suitable neighborhood of a fixed point $x_0$ where $u$ is differentiable. We also provide an estimate of the difference between $F^{**}(u)$ and the value $F^{**}(\tilde{u})$ in the modified function. In the last step of the proof we use the Vitali Covering

Theorem to obtain the desired function $v$. The proof is inspired by by a technique due to Cellina [22] in the framework of non-convex problems. The main difficulty we have to face here is that we do not assume, as in [22] and subsequent papers [45, 46, 47], that $(\nabla u(x_0), f^{**}(\nabla u(x_0)))$ is in the interior of a $N$-dimensional face of the epigraph of $f^{**}$.

In the whole proof we assume, without restriction, that $\int_\Omega f^{**}(\nabla u)\,dx < +\infty$.

*Step 1. Local construction in the case $k = N$.*

We remark that in this case we can argue as in [46] and take advantage of the construction performed there. We repeat the construction here for the sake of clarity and, moreover, to later underline the difficulties of the case $k < N$.

We fix a point $x_0 \in \Omega$ such that $u$ is differentiable at $x_0$ and we denote $\xi = \nabla u(x_0)$. Let $\xi_i$, $\alpha_i$ $1 \leq k \leq N$, the vectors and the coefficients involved in (HF) (ii). We define

$$\frac{1}{\alpha} = \max\{|x| : \max_{1,\ldots,N+1}(\xi_i - \xi) \cdot x \leq 1\} \qquad \text{and} \qquad \beta = \max_{i=1,\ldots N+1} |\xi_i - \xi|.$$

and we notice that assumption (HF) (ii) implies $\alpha > 0$ and that

$$\alpha |x - x_0| \leq \max_{i=1,\ldots,N+1}(\xi_i - \xi) \cdot (x - x_0) \leq \beta |x - x_0| \quad \forall x \in \mathbb{R}^N.$$

Now, for every $s \in \mathbb{R}$, $s > 0$, we define

$$\tilde{\omega}_{x_0,s}(x) := u(x_0) + \xi \cdot (x - x_0)$$
$$+ \max_{i=1,\ldots,N+1}(\xi_i - \xi) \cdot (x - x_0) - \frac{1}{2}s$$

and

$$\omega_{x_0,s}(x) := \tilde{\omega}_{x_0,s}(x)\, I_{B(x_0,s)}(x)$$

where

$$I_{B_s(x_0)}(x) = \begin{cases} 1 & \text{if } x \in \overline{B(x_0, s)} \\ +\infty & \text{otherwise.} \end{cases}$$

We claim that there exists $\bar{s}$ such that for every $0 < s < \bar{s}_{x_0} = \min\{\bar{s}, 2\alpha\bar{s}\}$, we have

$$u(x) \leq \omega_{x_0,s}(x) \qquad \text{on} \quad \partial B\left(x_0, \frac{s}{\alpha}\right). \tag{3.1}$$

and

$$\omega_{x_0,s}(x) \leq u(x) \qquad \text{on} \quad B\left(x_0, \frac{s}{(2\beta + \alpha)}\right). \tag{3.2}$$

To prove (3.1) and (3.2) we first recall that the differentiability property of $u$ in $x_0$ implies that there exists $\bar{s} > 0$ such that $B(x_0, \bar{s}) \subset \Omega$ and, for every $x \in B(x_0, \bar{s})$, we have

$$u(x_0) + \xi \cdot (x - x_0) - \frac{\alpha}{2}|x - x_0| \qquad (3.3)$$
$$< u(x) < u(x_0) + \xi \cdot (x - x_0) + \frac{\alpha}{2}|x - x_0|.$$

It follows, for every $x \in \partial B\left(x_0, \frac{s}{\alpha}\right)$,

$$u(x) \leq u(x_0) + \xi \cdot (x - x_0) + \frac{\alpha}{2}|x - x_0|$$
$$\leq u(x_0) + \xi \cdot (x - x_0) + \frac{1}{2} \max_{1,\ldots,N+1} (\xi_i - \xi) \cdot (x - x_0)$$
$$= u(x_0) + \xi \cdot (x - x_0) + \max_{1,\ldots,N+1} (\xi_i - \xi) \cdot (x - x_0)$$
$$- \frac{1}{2} \max_{1,\ldots,N+1} (\xi_i - \xi) \cdot (x - x_0)$$
$$\leq u(x_0) + \xi \cdot (x - x_0) + \max_{1,\ldots,N+1} (\xi_i - \xi) \cdot (x - x_0)$$
$$- \frac{\alpha}{2}|x - x_0|.$$

and the proof of (3.1) is completed observing that if $x \in \partial B\left(x_0, \frac{s}{\alpha}\right)$ the last expression coincides with $\omega_{x_0,,s}(x)$.

To prove (3.2) we consider $x \in B\left(x_0, \frac{s}{2\beta+\alpha}\right)$ and we compute

$$u(x_0) + \xi \cdot (x - x_0) + \max_{1,\ldots,N+1} (\xi_i - \xi) \cdot (x - x_0) - \frac{s}{2}$$
$$\leq u(x_0) + \xi \cdot (x - x_0) + \beta|x - x_0| - \frac{s}{2}$$
$$\leq u(x_0) + \xi \cdot (x - x_0) + \beta|x - x_0| - \frac{1}{2}(2\beta + \alpha)|x - x_0|$$
$$= u(x_0) + \xi \cdot (x - x_0) - \frac{\alpha}{2}|x - x_0| \leq u(x).$$

Inequality (3.1) implies that the function

$$\tilde{u}_{x_0,s} := \min\{u, \omega_{x_0,s}(x)\}$$

coincides with $u$ in $\Omega \setminus B\left(x_0, \frac{s}{\alpha}\right)$, hence also on $\partial\Omega$.

We denote by $E_{x_0,s}$ the set $\{x \in B\left(x_0 \frac{s}{\alpha}\right) : \tilde{u}_{x_0,s} = \omega_{x_0,s}(x)\}$ and we notice that, (3.1) and (3.2)

$$B\left(x_0, \frac{s}{2\beta+\alpha}\right) \subset E_{x_0,s} \subset B\left(x_0, \frac{s}{\alpha}\right). \tag{3.4}$$

Now we want to estimate $F^{**}(\tilde{u}_{x_0,s})$. A key point is the property $(b)$ in Lemma 2.4. We observe that, since

$\tilde{u}_{x_0,s} = u$, and hence also $\nabla \tilde{u}_{x_0,s} = \nabla u$, on $\Omega \setminus E_{x_0,s}$ it is sufficient to compare

$$\int_{E_{x_0,s}} f^{**}(\nabla u)\, dx \quad \text{and} \quad \int_{E_{x_0,s}} f^{**}(\nabla \tilde{u}_{x_0,s})\, dx.$$

The convexity of $f^{**}$, assumption (HF) (ii) and Lemma 2.4 imply that we can choose a selection $p(\cdot) \in \partial f^{**}(\cdot)$ such that $p(\nabla \tilde{u}_{x_0,s}(x)) = a$ for a.e. $x \in \Omega$. It follows

$$\int_{E_{x_0,s}} f^{**}(\nabla u)\, dx$$

$$\geq \int_{E_{x_0,s}} f^{**}(\nabla \tilde{u}_{x_0,s})\, dx + \int_{E_{x_0,s}} a(\nabla u - \nabla \tilde{u}_{x_0,s})\, dx$$

$$= \int_{E_{x_0,s}} f^{**}(\nabla \tilde{u}_{x_0,s})\, dx + \int_{B(x_0,\frac{s}{\alpha})} a(\nabla u - \nabla \tilde{u}_{x_0,s})\, dx$$

$$= \int_{E_{x_0,s}} f^{**}(\nabla \tilde{u}_{x_0,s})\, dx.$$

We notice that the last two equalities follow from the fact that $\nabla u = \nabla \tilde{u}_{x_0,s}$ a.e on $B(x_0, \frac{s}{\alpha}) \setminus E_{x_0,s}$ and observing that $\mathrm{div}\, a = 0$.
Finally we also underline that assumption (HF) (ii) and (iii) imply

$$\int_{E_{x_0,s}} f^{**}(\nabla u)\, dx \geq \int_{E_{x_0,s}} f(\nabla \tilde{u}_{x_0,s})\, dx.$$

*Step 2. Construction of a suitable piecewise affine function in the case $k < N$.*
As in the previous step, we fix a point $x_0 \in \Omega$ such that $u$ is differentiable at $x_0$, we denote $\xi = \nabla u(x_0)$ and we fix $R$ such that $R > |\xi| + 1$. Let $\xi_i$, $\alpha_i$, $\zeta_j$ be the vectors and the coefficients that satisfy assumption (HF) (ii) and (HF) (iii) for $\xi$ and $R$.
To begin with we need to fix some parameters that are needed in the construction. As we will see, some of them will depend on the choice of $x_0$. In order to keep the notation light, we will not always make explicit the dependence on $x_0$. We fix $\varepsilon > 0$. We assumed $f^{**}(\nabla u) \in L^1(\Omega)$ and so there exists $\delta > 0$ such that if $\Omega' \subset \Omega$ and $\lambda(\Omega') < \delta$ implies

$$\int_{\Omega'} |f^{**}(\nabla u)|\, dx < \frac{\varepsilon}{3\lambda(\Omega)}\,. \tag{3.5}$$

We still denote by $a$ the vector in $\mathbb{R}^N$ given by Lemma 2.4 (b) such that $a \in \partial f^{**}(\xi_i)$ for every $i = 1, \ldots, k+1$ and we choose $\gamma\, \eta \in \mathbb{R}$ such that

$$0 < \gamma < \min\left\{\frac{\varepsilon}{3^{N+2}(N-k)|a|\lambda(\Omega)}, 1\right\} \tag{3.6}$$

$$\eta = \frac{\delta}{\lambda(\Omega)}. \tag{3.7}$$

We set
$$M_{x_0} := \max_{i,j}\{|f(\xi_i)|, |f(\zeta_j)|, 1\} < +\infty \tag{3.8}$$

and we fix the real number $S_\eta$ such that

$$S_\eta \geq \max\left\{\frac{1}{1-\left(1-\frac{\varepsilon}{M_{x_0}\lambda(\Omega)3^{N+1}}\right)^{\frac{1}{N-k}}}, \frac{1}{1-(1-\frac{\eta}{3^N})^{\frac{1}{N-k}}}, \frac{3}{2}\right\}. \tag{3.9}$$

Without losing generality, we assume that

$$\xi_i \in \mathrm{Span}\{e_1, \ldots, e_k\} =: V \qquad \forall i = 1, \ldots, k+1$$

and we denote
$$V^\perp := \mathrm{Span}\{e_{k+1}, \ldots, e_N\}.$$

Given a point $x := (x_1, \ldots, x_N) \in \mathbb{R}^N$ we denote by $x' = (x_1, \ldots, x_k)$ its projection on $V$ and by $x'' = (x_{k+1}, \ldots, x_N)$ its projection on $V^\perp$. We consider the function $w_k : V \to \mathbb{R}$ defined by

$$w_k(x') := \max_{1,\ldots,k+1} (\xi_i - \xi) \cdot x'.$$

and we will also use
$$|x''|_{\infty,V^\perp} = \max_{j=k+1,\ldots,N} |x_j|$$

$$B := \{x' \in V : w_k(x') < \gamma\}$$
$$Q := \{x = (x', x'') \in \mathbb{R}^N : w_k(x') < \gamma \text{ and } |x''|_{\infty, V^\perp} < S_\eta\},$$
$$\tilde{Q} := \{x = (x', x'') \in \mathbb{R}^N : w_k(x') < \gamma \text{ and } |x''|_{\infty, V^\perp} < S_\eta - 1\}.$$

In the following, we will use that

$$\lambda(Q) = \lambda_k(B) S_\eta^{N-k} \tag{3.10}$$
$$\lambda(\tilde{Q}) = \lambda_k(B)(S_\eta - 1)^{N-k} \tag{3.11}$$
$$\lambda_{N-1}(\partial \tilde{Q} \setminus \partial Q) = 2(N-k)\lambda_k(B) S_\eta^{N-k-1} \tag{3.12}$$

We introduce the function $w : \mathbb{R}^N \to \mathbb{R}$ defined by

$$w(x) := \max\{w_k(x'), \frac{\gamma}{S_\eta}|x''|_{\infty, V^\perp}\}$$

Using the positive $1-$homogeneity of $w_k$ and $|\cdot|_{\infty, V^\perp}$ and the definition of $Q$, it follows that the function $w$ is 1-positively homogeneous, and satisfies the following properties:

$$\{x \in \mathbb{R}^N : w(x) < \gamma\} = Q;$$

$$\text{there exist} \quad 0 < \alpha < \beta \quad \text{such that} \quad \alpha|x| \leq w(x) \leq \beta|x| \quad \forall x \in \mathbb{R}^N \tag{3.13}$$

where the last property can be obtained in a similar way to the analogous one of Step 1.

We define the function $\tilde{w} : \mathbb{R}^N \to \mathbb{R}$ defined by

$$\tilde{w}(x) := \max\{w_k(x'), \gamma(|x''|_{\infty, V^\perp} - (S_\eta - 1))\}$$

and we prove that

$$\tilde{w}(x) \leq w(x) \quad \text{in} \quad Q, \tag{3.14}$$
$$\tilde{w}(x) = w(x) \quad \text{on} \quad \partial Q \tag{3.15}$$
$$\tilde{w}(x) = w_k(x) \quad \text{in} \quad \tilde{Q}. \tag{3.16}$$

In fact, since $w_k \geq 0$ and $|x''|_{\infty, V^\perp} \leq S_\eta - 1$ in $\tilde{Q}$, then (3.16) follows. To prove (3.14) and (3.15), it is sufficient to note that if $|x''| \leq S_\eta$

$$\gamma(|x''|_{\infty, V^\perp} - (S_\eta - 1)) \leq \frac{\gamma}{S_\eta}|x''|_{\infty, V^\perp}$$

where the equality holds if (and only if) $|x''|_{\infty,V^\perp} = S_\eta$.

Now let $A$ be the set

$$A := \{x \in Q : \tilde{w}(x) = \gamma(|x''|_{\infty,V^\perp} - (S_\eta - 1))\}$$

and notice that, by the definition of $\tilde{w}$ and (3.16), it follows that

$$A \subset Q \setminus \tilde{Q} \quad \text{and} \quad |\nabla \tilde{w}(x)| = \gamma < 1 \text{ a.e. in } A.$$

Denoting by $\zeta_j$ the vectors chosen at the beginning of this step, since $B(\xi, 1) \subset B(0, R+1) \subset \overline{co}\,(\cup_{j=1}^{N+1}\{\zeta_j\})$, it follows that $\nabla \tilde{w}(x) \in B(0,1) \subset (\cup_{j=1}^{N+1}\{\zeta_j - \xi\})$ and hence we can use the construction performed by Cellina in [22] and then subsequently refined by many authors (see [46] and references therein) to obtain a function $\overline{w}_A : A \to \mathbb{R}$ such that

$$\tilde{w}(x) - \frac{1}{2} \leq \overline{w}_A(x) \leq \tilde{w}(x) \quad \text{for every } x \in A, \tag{3.17}$$

$$\overline{w}_A(x) = \tilde{w}(x) \quad \text{for every } x \in \partial A$$

and, observing that $0 \in \text{int}\,\overline{co}\{\zeta_j - \xi\}$

$$\nabla \overline{w}_A(x) = \zeta_j - \nabla u(x_0) \quad \text{a.e. in } A.$$

Now we consider the function $\overline{w} : \mathbb{R}^n \to \mathbb{R}$ defined by

$$\overline{w}(x) := \begin{cases} \tilde{w}(x) & \text{if } x \in \mathbb{R}^N \setminus A \\ \overline{w}_A(x) & \text{if } x \in A \,. \end{cases}$$

*Step 3. Local construction in the case $K < N$.*

In this step we use the function $\overline{w}$ defined in Step 2 to obtain the desired local modification of the function $u$ in a suitable neighborhood of the point $x_0$ fixed at the beginning of Step 2. First of all, we recall that in Step 1 we already remarked that, by the differentiability of $u$ at $x_0$, there exist $\overline{r} \in \mathbb{R}$ such that $B(x_0, \overline{r}) \subset \Omega$ and inequalities (3.3) holds for every $x \in B(x_0, \overline{r})$.

We consider the set

$$Q_{x_0,s} = \{x \in \mathbb{R}^n : w(x - x_0) < s\gamma\} = x_0 + sQ$$

and we have that there exists $\overline{s}_{x_0}$ such that, for every $0 < s < \overline{s}_{x_0}$, $Q_{x_0,s} \subset B(x_0, \overline{r})$. For every $s \in ]0, \overline{s}_{x_0}[$ it also holds:

$$u(x) \leq u(x_0) + \xi \cdot (x - x_0) + w(x - x_0) - \frac{s}{2}\gamma \quad \forall x \in \partial Q_{x_0,s} \tag{3.18}$$

and
$$u(x_0) + \xi \cdot (x - x_0) + w(x - x_0) - \frac{s}{2}\gamma \leq u(x) \quad \forall x \in Q_{x_0, \frac{s}{3}}. \qquad (3.19)$$

In fact, using the second inequality in (3.3) and (3.13), we get

$$\begin{aligned} u(x) &\leq u(x_0) + \xi \cdot (x - x_0) + \frac{\alpha}{2}|x - x_0| \\ &\leq u(x_0) + \xi \cdot (x - x_0) + \frac{1}{2}w(x - x_0) \\ &= u(x_0) + \xi \cdot (x - x_0) + w(x - x_0) - \frac{1}{2}w(x - x_0) \end{aligned}$$

and recalling that on $\partial Q_{x_0, s}$ we have $w(x - x_0) = s\gamma$ we obtain (3.18). To deduce (3.19) we notice that in $Q_{x_0, \frac{s}{3}}$ it holds $3w(x - x_0) \leq s\gamma$; we use the first inequality in (3.3) and (3.13), so that

$$\begin{aligned} & u(x_0) + \xi \cdot (x - x_0) + w(x - x_0) - \frac{s}{2}\gamma \\ &\leq u(x_0) + \xi \cdot (x - x_0) + w(x - x_0) - \frac{3}{2}w(x - x_0) \\ &\leq u(x_0) + \xi \cdot (x - x_0) - \frac{\alpha}{2}|x - x_0| \leq u(x). \end{aligned}$$

Now we define
$$\overline{w}_s(x) = s\overline{w}\left(\frac{x}{s}\right)$$

and we note that $\overline{w}_s(x) = s\gamma$ for every $x \in s\partial Q$, and $\overline{w}_s \leq w(x)$ for every $x \in sQ$. Hence, from (3.14) and (3.15), we also have that

$$u(x) \leq u(x_0) + \xi \cdot (x - x_0) + \overline{w}_s(x - x_0) - \frac{s}{2}\gamma \quad \forall x \in \partial Q_{x_0, s} \qquad (3.20)$$

and

$$u(x_0) + \xi \cdot (x - x_0) + \overline{w}_s(x - x_0) - \frac{s}{2}\gamma \leq u(x) \quad \forall x \in Q_{x_0, \frac{s}{3}}. \qquad (3.21)$$

We define
$$\tilde{\omega}_{x_0, s}(x) = u(x_0) + \xi \cdot (x - x_0) + \overline{w}_s(x - x_0) - \frac{s}{2}\gamma$$

and
$$\omega_{x_0, s} = \tilde{\omega}_{x_0, s} I_{Q_{x_0, s}}$$

where
$$I_{Q_{x_0,s}}(x) = \begin{cases} 1 & \text{in } Q_{x_0,s}, \\ +\infty & \text{otherwise}. \end{cases}$$

Analogously to Step 1, we consider the function $\tilde{u}_{x_0,s} : \Omega \to \mathbb{R}$ defined by
$$\tilde{u}_{x_0,s}(x) := \min\{u(x), \omega_{x_0,s}\}$$

and we denote by $E_{x_0,s}$ the set $\{x \in Q_{x_0,s} : \tilde{u}_{x_0,s} = \omega_{x_0,s}(x)\}$. By (3.20) and (3.21), it follows that
$$Q_{x_0, \frac{s}{3}} \subset E_{x_0,s} \subset Q_{x_0,s} \tag{3.22}$$

We observe that
$$\nabla \tilde{u}_{x_0,s}(x) = \nabla u(x) \quad \text{for a.e. } x \in Q_{x_0,s} \setminus E_{x_0,s}, \tag{3.23}$$

$$\nabla \tilde{u}_{x_0,s}(x) \in \{\xi_i, i = 1, \ldots, k+1\} \cup \{\zeta_j, j = 1, \ldots, N+1\} \tag{3.24}$$
$$\text{for a.e. } x \in \tilde{Q}_{x_0,s} \cap E_{x_0,s},$$

$$f(\nabla \tilde{u}_{x_0,s}(x)) = f^{**}(\nabla \tilde{u}_{x_0,s}(x)) \text{ a.e. } x \in E_{x_0,s}. \tag{3.25}$$

Analogously to Step 1, we compare
$$\int_{E_{x_0,s}} f^{**}(\nabla u) \, dx \quad \text{and} \quad \int_{E_{x_0,s}} f(\nabla \tilde{u}_{x_0,s}) \, dx.$$

We start by considering
$$\int_{E_{x_0,s}} f^{**}(\nabla u) \, dx \tag{3.26}$$
$$= \int_{E_{x_0,s} \cap (Q_{x_0,s} \setminus \tilde{Q}_{x_0,s})} f^{**}(\nabla u) \, dx - \int_{E_{x_0,s} \cap (Q_{x_0,s} \setminus \tilde{Q}_{x_0,s})} f^{**}(\nabla \tilde{u}_{x_0,s}) \, dx$$
$$+ \int_{E_{x_0,s} \cap (Q_{x_0,s} \setminus \tilde{Q}_{x_0,s})} f^{**}(\nabla \tilde{u}_{x_0,s}) \, dx + \int_{E_{x_0,s} \cap \tilde{Q}_{x_0,s}} f^{**}(\nabla u) \, dx$$

By the convexity of $f^{**}$ and Lemma 2.4 (b), we can choose a selection $p(\cdot) \in \partial f^{**}(\cdot)$ such that $p(\xi_i) = a$ for every $i = 1, \ldots, k+1$, so that also using (3.24) we get
$$\int_{E_{x_0,s} \cap \tilde{Q}_{x_0,s}} f^{**}(\nabla u) \, dx \tag{3.27}$$

$$\geq \int_{E_{x_0,s} \cap \tilde{Q}_{x_0,s}} f^{**}(\nabla \tilde{u}_{x_0,s})\, dx + \int_{E_{x_0,s} \cap \tilde{Q}_{x_0,s}} a \cdot (\nabla u - \nabla \tilde{u}_{x_0,s})\, dx$$
$$= \int_{E_{x_0,s} \cap \tilde{Q}_{x_0,s}} f(\nabla \tilde{u}_{x_0,s})\, dx + \int_{\tilde{Q}_{x_0,s}} a \cdot (\nabla u - \nabla \tilde{u}_{x_0,s})\, dx$$

where, in the last equality, we applied (3.23). The fact that $a$ is constant on $\tilde{Q}_{x_0,s}$ implies that

$$\int_{\tilde{Q}_{x_0,s}} a \cdot (\nabla u - \nabla \tilde{u}_{x_0,s})\, dx = \int_{\partial \tilde{Q}_{x_0,s}} (u - \tilde{u}_{x_0,s}) a \cdot \nu_x dH^{N-1} \qquad (3.28)$$

where $\nu_x$ denotes the external normal to $\tilde{Q}_{x_0,s}$. By construction and (3.20), recalling that $B = \{x' \in V : w_k(x') < \gamma\}$, it follows that

$$\{x \in \partial \tilde{Q}_{x_0,s} : u(x) \neq \tilde{u}_{x_0,s}(x)\} \subset \partial \tilde{Q}_{x_0,s} \setminus \partial Q_{x_0,s}$$
$$= \{x \in \mathbb{R}^n : x' \in sB \text{ and } |x''|_{\infty,V^\perp} = s(S_\eta - 1)\}$$

thus, recalling (3.12),

$$\lambda_{N-1}(\{x \in \partial \tilde{Q}_{x_0,s} : u(x) \neq \tilde{u}_{x_0,s}(x)\})$$
$$\leq \lambda_{N-1}(s(\partial \tilde{\Omega} \setminus \partial Q))$$
$$= 2(N-k)s^{N-1}\lambda_k(B)(S_\eta - 1)^{N-k-1}.$$

Furthermore, using (3.13), (3.17) and (3.3), we have, on $\partial \tilde{Q}_{x_0,s}$,

$$|\tilde{u}_{x_0,s}(x) - u(x)| \leq \frac{3}{2} s\gamma$$

so that, returning to (3.28), we have the following

$$\left| \int_{\tilde{Q}_{x_0,s}} a \cdot (\nabla u - \nabla \tilde{u}_{x_0,s}) dx \right|$$
$$\leq 3\gamma(N-k)|a|s^N \lambda_k(B)(S_\eta - 1)^{N-k-1}$$
$$= 3^{N+1}\gamma(N-k)|a|\frac{s^N \lambda_k(B)(S_\eta - 1)^{N-k-1}}{3^N}$$
$$= 3^{N+1}\gamma(N-k)|a|\lambda(Q_{x_0,\frac{s}{3}})$$
$$\leq 3^{N+1}\gamma(N-k)|a|\lambda(E_{x_0,s}).$$

Recalling (3.6) we get

$$\left| \int_{E_{x_0,s} \cap \tilde{Q}_{x_0,s}} a \cdot (\nabla u - \nabla \tilde{u}_{x_0,s}) \, dx \right| \leq \frac{\varepsilon}{3} \frac{\lambda(E_{x_0,s})}{\lambda(\Omega)} \, . \tag{3.29}$$

Now we want to estimate the measure of the set of integration in the expression

$$\int_{E_{x_0,s} \cap (Q_{x_0,s} \setminus \tilde{Q}_{x_0,s})} f^{**}(\nabla u) \, dx \, .$$

By the choice of $S_\eta$ in (3.9) and by (3.10) and (3.11), we obtain

$$\begin{aligned} \lambda(E_{x_0,s} \cap (Q_{x_0,s} \setminus \tilde{Q}_{x_0,s})) &\leq \lambda(Q_{x_0,s} \setminus \tilde{Q}_{x_0,s}) \\ &= s^N \lambda(Q \setminus \tilde{Q}) \\ &= s^N \lambda_k(B) \left( S_\eta^{N-k} - (S_\eta - 1)^{N-k} \right) \\ &\leq 3^N \left(\frac{s}{3}\right)^N \lambda_K(B) \left( S_\eta^{N-k} - (S_\eta - 1)^{N-k} \right) \\ &\leq \eta \lambda(Q_{x_0, \frac{s}{3}}) \leq \delta \frac{\lambda(E_{x_0,s})}{\lambda(\Omega)}. \end{aligned} \tag{3.30}$$

It remains to estimate

$$\int_{E_{x_0,s} \cap (Q_{x_0,s} \setminus \tilde{Q}_{x_0,s})} f^{**}(\nabla \tilde{u}_{x_0,s}) \, dx.$$

Using (3.8) and (3.25) we obtain

$$|f^{**}(\nabla \tilde{u}_{x_0,s}(x))| = |f(\nabla \tilde{u}_{x_0,s}(x))| \leq M_{x_0}$$

for a.e. $x \in E_{x_0,s} \cap (Q_{x_0,s} \setminus \tilde{Q}_{x_0,s})$, so that, using once again (3.9), (3.10) and (3.11), we obtain

$$\begin{aligned} \left| \int_{E_{x_0,s} \cap (Q_{x_0,s} \setminus \tilde{Q}_{x_0,s})} f^{**}(\nabla \tilde{u}_{x_0,s}) \, dx \right| \\ \leq M_{x_0} \lambda(Q_{x_0,s} \setminus \tilde{Q}_{x_0,s}) = M_{x_0} s^N \lambda_K(B) \left( S_\eta^{N-k} - (S_\eta - 1)^{N-k} \right) \\ \leq \frac{\varepsilon}{3} \frac{\lambda(E_{x_0,s})}{\lambda(\Omega)} \end{aligned} \tag{3.31}$$

Collecting (3.26), (3.27), (3.29), (3.30), (3.31) and recalling that $f^{**}(\nabla \tilde{u}_{x_0,s}(x)) = f((\nabla \tilde{u}_{x_0,s}(x))$ for a.e. $x \in E_{x_0,s}$ we obtain that

$$\int_{E_{x_0,s}} f(\nabla \tilde{u}_{x_0,s}) \, dx = \int_{E_{x_0,s}} f^{**}(\nabla \tilde{u}_{x_0,s}) \, dx \qquad (3.32)$$

$$\leq \int_{E_{x_0,s}} f^{**}(\nabla u) \, dx$$

$$+ \left| \int_{E_{x_0,s} \cap (Q_{x_0,s} \setminus \tilde{Q}_{x_0,s})} f^{**}(\nabla u) \, dx \right| + 2 \frac{\varepsilon}{3} \frac{\lambda(E_{x_0,s})}{\lambda(\Omega)}.$$

*Step 4. Construction of the function $v$.*
We consider the bounded measurable set

$$\tilde{\Omega} := \{x \in \Omega : u \text{ is differentiable at } x \text{ and } f(\nabla u(x)) \neq f^{**}(\nabla u(x))\}.$$

For every $x \in \tilde{\Omega}$ we can consider the family of sets $E_{x,s}$, $0 < s < \bar{s}_x$ determined in Step 1 for the case where $\nabla u(x)$ satisfies assumption (HF) (ii) with $k = N$ or in Step 3 for the case $k < N$.
We prove, as in Claim 2.5 of [46], that for every $E_{x,s}$ there exists a closed set $F_{x,s}$ such that

$$\lambda(E_{x,s}) = \lambda(F_{x,s}) \qquad (3.33)$$

and

$$E_{x,s} \subseteq F_{x,s} \subseteq B\left(x, \frac{s}{\alpha}\right) \qquad \text{if} \quad k = N \qquad (3.34)$$

$$E_{x,s} \subseteq F_{x,s} \subseteq Q_{x,s} \qquad \text{if} \quad k < N.$$

To this aim, we consider

$$G_{x,s} = \begin{cases} \{x \in B(x,s) : u(x) > \omega_{x,s}(x)\} & \text{if } k = N \\ \{x \in Q_{x,s} : u(x) > \omega_{x,s}(x)\} & \text{if } k < N \end{cases}$$

The assumption that $u$ is a.e. differentiable implies that, for a.e. $y \in G_{x,s}$, there exists $B(y, r_y)$ such that $B(y, s) \subseteq G_{x,s}$. Defining

$$F_{x,s} = \begin{cases} B(x,s) \setminus \cup_{y \in G_{x,s}} B(y,s) & \text{if } k = N \\ Q_{x,s} \setminus \cup_{x \in G_{x,s}} B(y,s) & \text{if } k < N \end{cases}$$

we completed the proof of (3.33) and (3.34) and recalling also (3.4) and (3.22) we have that the family $F_{x,s}$, for $x$ varying in $\tilde{Q}$, is a Vitali covering of $\tilde{Q}$, (see Theorem 3.1 in [43]). Then there exists an at most countable family $(F_{x_n,s_n}) = (F_n)_{n\in\mathbb{N}}$ of mutually disjoint elements of $\mathcal{F}$ such that

$$\lambda\left(\tilde{\Omega} \setminus \cup_{n\in\mathbb{N}} F_n\right) = \lambda\left(\tilde{\Omega} \setminus \cup_{n\in\mathbb{N}} E_n\right) = 0$$

where we denoted $E_n := E_{x_n,s_n}$.
Finally we can define

$$v(x) = \begin{cases} \tilde{u}_{x_n}(x) & \text{if } x \in F_n \\ u(x) & \text{if } x \notin \cup_n F_n \, . \end{cases}$$

Exploiting (3.32) we obtain

$$\begin{aligned}
\int_\Omega f(\nabla v)dx &= \int_{\Omega \setminus \cup_n F_n} f^{**}(\nabla u)dx + \int_{\cup_n F_n} f^{**}(\nabla v)dx \\
&= \int_{\Omega \setminus \cup_n F_n} f^{**}(\nabla u)dx + \sum_n \int_{F_n} f^{**}(\nabla v)dx \\
&\leq \int_{\Omega \setminus \cup_n F_n} f^{**}(\nabla u)dx + \sum_n \int_{F_n} f^{**}(\nabla u)dx \\
&\quad + \sum_n \left| \int_{E_n \cap (Q_{x_n,s_n} \setminus \tilde{Q}_{x_n,s_n})} f^{**}(\nabla v)\, dx \right| \\
&\quad + \frac{2}{3}\varepsilon \sum_n \frac{\lambda(E_{x_n,s_n})}{\lambda(\Omega)} \\
&\leq \int_\Omega f^{**}(\nabla u(x)dx + \varepsilon.
\end{aligned}$$

We notice that in the last inequality we used the estimate (3.30) on the measure of the sets $E_n \cap (Q_{x_n,s_n} \setminus \tilde{Q}_{x_n,s_n})$ and the choices of $\eta$ and $\delta$ in (3.7) and (3.5). To complete the proof it is sufficient to notice that, by construction, we have $v = \varphi$ on $\partial\Omega$ and that the last inequality, together with assumption (HF) (i) and Poincaré inequality, implies that $v \in \varphi + W_0^{1,1}(\Omega)$.

□

# 4 Non-occurrence of the Lavrentiev gap

The aim of this section is to apply the construction of Section 3 to obtain a result on the non-occurrence of the Lavrentiev gap. First of all, the next theorem is a consequence of Theorem 3.1 that shows that if $u \in \varphi + W_0^{1,\infty}(\Omega)$ then there exists a sequence $v_n \in \varphi + W_0^{1,\infty}(\Omega)$ such that approximates $u$ in energy,
Here we need to slightly modify the assumption (HF) requiring that the growth of the Lagrangian is superlinear. To be more precise, we formulate the following

(HF')  (i) There exist $\phi : \mathbb{R}^+ \to \mathbb{R}$ such that $\lim_{t \to +\infty} \frac{\phi(t)}{t} = +\infty$ and

$$f(\zeta) \geq \phi(|\zeta|) \qquad \text{for every } \zeta \in \mathbb{R}^N$$

(ii) For every $\xi \in \mathbb{R}^N$ either $f^{**}(\xi) = f(\xi)$ or there exist $k \in \{1, \ldots, N\}$, $\{\xi_i\}_{i=1\ldots,k+1} \subset \mathbb{R}^N$ and $\{\alpha_i\}_{i=1\ldots,k+1} \subset \mathbb{R}$ such that $\alpha_i$ is strictly positive for every $i = 1, \ldots, k+1$, $\sum_{i=1}^{k+1} \alpha_i = 1$, $\sum_{i=1}^{k+1} \alpha_i \xi_i = \xi$ and

$$\sum_{i=1}^{k+1} \alpha_i f(\xi_i) = f^{**}(\xi).$$

Moreover we also assume that $\dim \mathrm{Span}(\xi_1, ..., \xi_{k+1}) = k$.

(iii) For every $R > 0$ there exist $N+1$ vectors $\zeta_j \in \mathbb{R}^N$ such that $\overline{B}(0, R) \subset \overline{\mathrm{co}}(\cup_{j=1}^{N+1}\{\zeta_j\})$ and

$$f(\zeta_j) = f^{**}(\zeta_j) < +\infty \qquad \text{for every } j = 1, ..., N+1.$$

**Remark 4.1.** *By [33][Theorem 4.98], if $f$ is superlinear and lower semicontinuous, then it satisfies (HF') (ii).*

Now we can prove the following theorem.

**Theorem 4.2.** *Let $\Omega$ be an open bounded Lipschitz subset of $\mathbb{R}^N$ and assume that $f : \mathbb{R}^N \to \mathbb{R} \cup \{+\infty\}$ satisfies assumption* (HF'). *Then for every $u \in W^{1,\infty}(\Omega)$ there exists $u_n \in u + W_0^{1,\infty}(\Omega)$ such that*

$$u_n \overset{*}{\rightharpoonup} u \quad \text{in} \quad W^{1,\infty}(\Omega)$$

*and*

$$\lim_n \int_\Omega f(\nabla u_n) dx = \int_\Omega f^{**}(\nabla u) dx.$$

**Proof.** First of all we notice that, by Lemma 2.4 $(a)$, it follows that $\int_\Omega f^{**}(\nabla u)\, dx$ is finite. Moreover $u$ is a.e. differentiable in $\Omega$, so that we can repeat the same construction that we made in the proof of Theorem 3.1.

We fix $u \in W^{1,\infty}(\Omega)$ and $R = \|\nabla u\|_\infty + 1$. We can consider the set

$$\Omega' = \{x \in \Omega : f^{**}(\nabla u(x)) \neq f(\nabla u(x)) \text{ and } u \text{ is differentiable at } x\}$$

For every $x \in \Omega'$ there exist $k(x) \in \{1,\ldots,N\}$, $\xi_i(x)$, $\alpha_i(x)$ satisfying (HF) (ii). Moreover (HF) (iii) states that there exist $\zeta_j$, $j = 1,\ldots,N+1$, such that

$$\overline{B}(0,R) \subset \overline{\operatorname{co}}(\cup_{j=1}^{N+1}\{\zeta_j\})$$

and we notice that the vectors $\zeta_j$ do not depend on the choice of $x$.

Using Theorem 3.1 we deduce that, for every $n \in \mathbb{N}$, there exists $v_n \in u+W_0^{1,1}(\Omega)$ such that

$$\int_\Omega f(\nabla v_n)dx \leq \int_\Omega f^{**}(\nabla u)dx + \frac{1}{n}$$

We claim that $v_n \in u+W_0^{1,\infty}(\Omega)$ and that there exists $K \in \mathbb{R}$ such that $\|\nabla v_n\|_\infty \leq K$ for every $n \in \mathbb{N}$. To prove this claim it is sufficient to recall that in the proof of Theorem 3.1 we obtain that $\nabla v_n(x) \in \{\cup_{i=1,\ldots,k(x)+1}\xi_i(x), \cup_{j=1,\ldots,N+1}\zeta_j\}$, for a.e. $x \in \Omega'$, and then it is sufficient to show that there exists $M > 0$ such that for every $x \in \Omega'$

$$|\xi_i(x)| \leq M \qquad \text{for every } i = 1,\ldots,k(x)+1.$$

To this aim we fix $x \in \Omega'$ and, in order to keep the notation light, we denote $\xi = \nabla u(x)$ and we drop the dependence on $x$ in the $\xi_i$ and $k$. By Lemma 2.4 (b) it follows that there exists $a \in \partial f^{**}(\xi)$ such that

$$f^{**}(\xi_i) = f^{**}(\xi) + a \cdot (\xi_i - \xi) \qquad \text{for every } i = 1,\ldots k+1. \qquad (4.1)$$

We define

$$g(\zeta) := \max_{b \in \partial f^{**}(\eta), \eta \in \overline{B(0,\|\nabla u\|_\infty)}} f^{**}(\eta) + b \cdot (\zeta - \eta)$$

and we observe that, by (4.1) and the convexity of $f^{**}$, it turns out that,

$$g(\xi_i) = f^{**}(\xi_i) \qquad \text{for every } i = 1,\ldots k+1. \qquad (4.2)$$

The inequality

$$g(\zeta) \leq \max_{\eta \in \overline{B(0,\|\nabla u\|_\infty)}} f^{**}(\eta) + \max_{b \in \partial f^{**}(\eta), \eta \in \overline{B(0,\|\nabla u\|_\infty)}} |b|(|\zeta| + \|\nabla u\|_\infty)$$

shows that the function $g$ grows at most linearly and hence, by the superlinearity of $f$ and $f^{**}$ there exists $M > 0$ such that if $|\zeta| > M$ then $g(\zeta) < f^{**}(\zeta)$. Thus, by (4.2), the claim is proved and we have that $\|\nabla v_n\|_\infty \leq K := \max\{M, |\zeta_j|, j = 1, \ldots, N+1\}$.

A standard argument implies the existence of a subsequence that we still denote by $v_n$ such that $v_n \stackrel{*}{\rightharpoonup} v$ and, by the weak*-lower semicontinuity of $F^{**}$, we obtain

$$\int_\Omega f^{**}(\nabla u) dx \leq \liminf_n \int_\Omega f^{**}(\nabla v_n) dx = \liminf_n \int_\Omega f(\nabla v_n) dx$$
$$\leq \limsup_n \int_\Omega f(\nabla v_n) dx \leq \int_\Omega f^{**}(\nabla u) dx$$

□

Now we can apply the previous result to prove the non-occurrence of the so-called Lavrentiev gap.

**Theorem 4.3.** *Let $\Omega$ be an open, bounded, Lipschitz subset of $\mathbb{R}^N$. Let $\varphi$ be in $W^{1,\infty}(\Omega)$ and assume that $f : \mathbb{R}^N \to \mathbb{R} \cup \{+\infty\}$ satisfies assumption (HF'). Then for every $u \in \varphi + W_0^{1,1}(\Omega)$ there exists $u_n \in \varphi + W_0^{1,\infty}(\Omega)$ such that*

$$u_n \rightharpoonup u \quad in \quad W^{1,1}(\Omega)$$

*and*

$$\lim_n \int_\Omega f(\nabla u_n) dx = \int_\Omega f^{**}(\nabla u) dx.$$

**Proof.** Is sufficient to consider the case where $\int_\Omega f^{**}(\nabla u)\, dx$. Theorem 1.1 in [10] and Theorem 5 in [12] state that for every $u \in \varphi + W_0^{1,1}(\Omega)$ there exists a sequence in $u_n \in \varphi + W^{1,\infty}(\Omega)$ such that $u_n$ strongly converge in $W^{1,1}(\Omega)$ to $u$ and

$$\lim_n \int_\Omega f^{**}(\nabla v_n) dx = \int_\Omega f^{**}(\nabla u) dx.$$

By Theorem 4.3, for every $v_n$ we can construct a sequence $\{v_n^h\}_{h \in \mathbb{N}}$ such that $v_n^h \in \varphi + W_0^{1,\infty}(\Omega)$ for every $h$, and

$$v_n^h \stackrel{*}{\rightharpoonup} v_n \quad in \quad W^{1,\infty}(\Omega),$$

$$\lim_n \int_\Omega f(\nabla v_n^h) dx = \int_\Omega f^{**}(\nabla v_n) dx.$$

Thus, via a diagonal argument, we can determine a sequence $\{u_n\}_{n\in\mathbb{N}}$ such that

$$\lim_n \int_\Omega f(\nabla u_n) dx = \int_\Omega f^{**}(\nabla u) dx$$

and then we also have that

$$\int_\Omega \phi(|\nabla u_n|)\, dx \leq \int_\Omega f(\nabla u_n) dx \leq \int_\Omega f^{**}(\nabla u) dx + C$$

for a suitable $C \in \mathbb{R}$. By de la Vallée Poussin Theorem we can conclude that, up to a subsequence, $u_n \rightharpoonup u$ in $W^{1,1}(\Omega)$. □

As a consequence, we can also prove the non-occurrence of the Lavrentiev phenomenon for the functional $F$. We will argue as in Theorem 30 in [6].

**Theorem 4.4.** *Let $\Omega$ be an open, bounded and Lipschitz subset of $\mathbb{R}^N$, let $f : \mathbb{R}^N \to \mathbb{R} \cup \{+\infty\}$ satisfy* (HF') *and $\varphi$ be in $W^{1,\infty}(\Omega)$. Then*

$$\inf_{\varphi+W_0^{1,1}(\Omega)} \int_\Omega f(\nabla u) dx = \inf_{\varphi+W_0^{1,\infty}(\Omega)} \int_\Omega f(\nabla u) dx\,.$$

**Proof.** It is sufficient to observe that

$$\inf_{\varphi+W_0^{1,\infty}(\Omega)} \int_\Omega f(\nabla u) dx = \inf_{\varphi+W_0^{1,1}(\Omega)} \int_\Omega f^{**}(\nabla u) dx$$

$$\leq \inf_{\varphi+W_0^{1,1}(\Omega)} \int_\Omega f(\nabla u) dx \leq \inf_{\varphi+W_0^{1,\infty}(\Omega)} \int_\Omega f(\nabla u) dx$$

where the first equality has been proved in Theorem 4.3. □

# References

[1] G. Alberti, Giovanni, P. Majer, *Gap phenomenon for some autonomous functionals*, J. Convex Anal., 1 (1) 1994, 31-45.

[2] G. Alberti, F. Serra Cassano, *Non-occurrence of gap for one-dimensional autonomous functionals*, Calculus of variations, homogenization and continuum mechanics (Marseille, 1993), 1–17. Ser. Adv. Math. Appl. Sci., (18) 1994, 1-17.

[3] A. Kh. Balci, L. Diening, M. Surnachev, *New examples on the Lavrentiev gap using fractals*, Calc. Var. & PDE, 5 (59) 2020, https://doi.org/10.1007/s00526-020-01818-1.

[4] A. Kh. Balci, L. Diening, M. Surnachev, *Scalar minimizers with maximal singular sets and lack of the Meyers property*, 2023, preprint (arXiv:2312.15772).

[5] P. Baroni, M. Colombo, G. Mingione, *Regularity for general functionals with double phase*, Calc. Var. Partial Diff. Equ. 2 (57) 2018, https://doi.org/10.1007/s00526-018-1332-z.

[6] T. Bertin, *Integral representations of lower semicontinuous envelopes and Lavrentiev phenomenon for non continuous Lagrangians*, Nonlinear Analysis: Real World Applications, (86) 2025, https://doi.org/10.1016/j.nonrwa.2025.104414.

[7] M. Borowski, I. Chlebicka, B. Miasojedow, *Absence of Lavrentiev's gap for anisotropic functionals* Nonlinear Anal., (246) 2024, https://doi.org/10.1016/j.na.2024.113584.

[8] M. Borowski, I. Chlebicka, F. De Filippis, and B. Miasojedow, *Absence and presence of Lavrentiev's phenomenon for double phase functionals upon every choice of exponents*, Calc. Var. Partial Differential Equations, 2 (63) 2024, https://doi.org/10.1007/s00526-023-02640-1.

[9] M. Borowski, P. Bousquet, I. Chlebicka, B. Lledos, B. Miasojedow. *Discarding Lavrentiev's Gap in Non-automonous and Non-Convex Variational Problems*, 2024 (to appear). hal-04814888.

[10] P. Bousquet, *Nonoccurence of the Lavrentiev gap for multidimensional autonomous problems*, Ann. Sc. Norm. Super. Pisa Cl. Sci. 5, 3 (24) 2023, 1611–1670.

[11] P. Bousquet, C. Mariconda, G. Treu, *On the Lavrentiev phenomenon for multiple integral scalar variational problems*, Journal of Functional Analysis, 9 (266) 2014, 5921-5954.

[12] P. Bousquet, C. Mariconda, G. Treu, *Non occurrence of the Lavrentiev gap for a class of nonautonomous functionals*, Ann. Mat. Pura Appl. (4), 52 (203) 2024, 2275–2317.

[13] P. Bousquet, C. Mariconda, G. Treu, *Some recent results on the Lavrentiev phenomenon*, Commun. Optimiz. Theory, 6 2024, 1–25.

[14] M. Bulícek, P. Gwiazda, J. Skrzeczkowski, *On a range of exponents for absence of Lavrentiev phenomenon for double phase functionals*, Arch. Ration. Mech. Anal., 1 (246) 2022, 209–240 .

[15] G. Buttazzo, *Semicontinuity, relaxation and integral representation in the calculus of variations*, Longman Scientific & Technical, (207) 1989.

[16] G. Buttazzo, G. Dal Maso, *Γ-limits of integral functionals*, J. Analyse Math., (37) 1980, 145–185.

[17] G. Buttazzo, G. Dal Maso, *Integral representation and relaxation of local functionals*, Nonlinear Analysis: Theory, Methods & Applications, 6 (9) 1985, 515-532.

[18] G. Buttazzo, G. Dal Maso, E. De Giorgi, *On the lower semicontinuity of certain integral functionals*, Atti della Accademia Nazionale dei Lincei. Classe di Scienze Fisiche, Matematiche e Naturali. Rendiconti Lincei. Matematica e Applicazioni, 5 (74) 1983, 274-282.

[19] G. Buttazzo, A. Leaci, *Relaxation Results for a Class of Variational Integrals*, Journal of Functional Analysis, 3 (61) 1985, pp 360-377.

[20] G. Buttazzo, V. J. Mizel, *Interpretation of the Lavrentiev phenomenon by relaxation*, J. Funct. Anal., 2 (110) 1992, 434–460.

[21] A. Cellina, *On minima of a functional of the gradient: necessary conditions*, Non Linear Analysis, 4 (20) 1993, pp 337-341.

[22] A. Cellina, *On minima of a functional of the gradient: sufficient conditions*, Non Linear Analysis, 4 (20) 1993, 343-347.

[23] B. Dacorogna, *Introduction to the calculus of variations*, Imperial College Press, London, 2015.

[24] G. Dal Maso, *Relaxation of Autonomous Functionals with Discontinuous Integrands*, nn. Univ. Ferrara - Sez. VII - 8c. Mat. (XXXIV) 1988, 21-47.

[25] R. De Arcangelis, *Some Remarks on the Identity between a Variational Integral and Its Relaxed Functional*, Ann. Univ. Ferrara - Sez. VII - Sc. Mat. (XXXV) 1989, 135-145.

[26] R. De Arcangelis, E. Zappale, *The Relaxation of Some Classes of Variational Integrals with Pointwise Continuous-Type Gradient Constraints*, Appl. Math. Optim., 3 (51) 2005, 251–277, .

[27] F. De Filippis, F. Leonetti, *No Lavrentiev gap for some double phase integrals*, Adv. Calc. Var. ; 1 (17) 2024, 165–194.

[28] F. De Filippis, F. Leonetti, G. Treu, *Nonoccurrence of Lavrentiev gap for a class of functionals with nonstandard growth*, Advances in Nonlinear Analysis 2024, 1 (13) 2024, https://doi.org/10.1515/anona-2024-0002.

[29] I. Ekeland, R. Témam, *Convex Analysis and Variational Problems*, 1999 (First edition 1976).

[30] A. Esposito, F. Leonetti, P. V. Petricca *Absence of Lavrentiev gap for non-autonomous functionals with (p, q)-growth*, Adv. Nonlinear Anal., 1 (8) 2019, 73–78.

[31] L. Esposito, F. Leonetti, G. Mingione, *Sharp regularity for functionals with $(p, q)$-growth*, J. Diff. Equ., 1 (204) 2004, 5-55.

[32] I. Fonseca, J. Malý, G. Mingione, *Scalar minimizers with fractal singular sets*, Arch. Ration. Mech. Anal., 2 (172) 2004, 295-307.

[33] I. Fonseca, G, Leoni, *Modern Methods in the Calculus of Variations: $L^p$ Spaces*, 2007.

[34] G. Friesecke, *A necessary and sufficient condition for nonattainment and formation of microstructure almost everywhere in scalar variational problems*, Proceedings of the Royal Society of Edinburgh: Section A Mathematics, 3 (124) 1994, 437-471.

[35] L. Koch, M. Ruf, M. Schäffner, *On the lavrentiev gap for convex, vectorial integral functionals*, J. Funct. Anal. 5 (288) 2025, https://doi.org/10.1016/j.jfa.2024.110793

[36] M. Lavrentiev, *Sur quelques problèmes du calcul des variations* Ann. Mat. Pura Appl., (4) 1926, 107–124.

[37] A.C. Heinricher, V. Mizel, *A new example of the Lavrentiev phenomenon*, SIAM J. Control Optim., 6 (26) 1988, 1490-1503.

[38] A. C. Heinricher,; V. C: Mizel, *The Lavrentiev Phenomenon for lnvariant Variational Problems*, Archive for rational mechanics and analysis, 1 (102) 1988, 57-93.

[39] F. Hüsseinov, *Relaxation and Nonoccurrence of the Lavrentiev Phenomenon for Nonconvex Problems* Acta Mathematica Sinica, English Series Jun., 6 (29) 2013, 1185–1198.

[40] B. Manià, *Sopra un esempio di Lavrentieff*, Boll. Un. Matem. Ital. (13) 1934, 147–153.

[41] P. Marcellini, C. Sbordone, *Semicontinuity problems in the Calculus of Variations*, Non linear Analysis, Theory, Methods & Applications 2 (4) 1980, 241-257.

[42] C. Mariconda, *Non-occurrence of gap for one-dimensional non-autonomous functionals*, Calc. Var. Partial Differential Equations 2 (62) 2023, https://doi.org/10.1007/s00526-022-02391-5.

[43] S. Saks *Theory of the Integral*, 2nd revised edn. With two additional notes by Stefan Banach, 1964.

[44] J. Serrin, *On the definition and properties of certain variational integrals*, Trans. Amer. Math. Soc. (101) 1961, 139–167.

[45] M. A. Sychev. *Characterization of homogeneous scalar variational problems solvable for all boundary data*, Proc. Roy. Soc. Edinburgh Sect., 3 (130) 2000, 611–631.

[46] M. A. Sychev, G. Treu, G. Colombo, *Stable minimizers of functionals of the gradient*, Proc. Roy. Soc. Edinburgh Sect., 5 (150) 2020, 2642–2655.

[47] S. Zagatti, *Minimization of functionals of the gradient by Baire's theorem*, SIAM J. Control Optim. 2 (38) 2000, 384–399.

[48] W. P. Ziemer. *Weakly differentiable functions. Sobolev spaces and functions of bounded variation. Graduate Texts in Mathematics*, New York: Springer, (120) 1989.

[49] V.V. Zhikov, *On Lavrentiev's phenomenon* Russ. J. Math. Phys. 2 (3) 1995, 249-269.